\documentclass[]{interact}

\usepackage{epstopdf}
\usepackage{subfigure}

\usepackage{natbib}
\usepackage{xcolor}
\bibpunct[, ]{(}{)}{;}{a}{}{,}
\renewcommand\bibfont{\fontsize{10}{12}\selectfont}

\theoremstyle{plain}

\theoremstyle{definition}

\theoremstyle{remark}

\def\cl {\nonumber \\}
\def\el {\nonumber }

\begin{document}


\title{Reduced Basis Model Order Reduction for Navier-Stokes equations in domains with walls of varying curvature}

\author{
\name{Martin~W. Hess\textsuperscript{a}\thanks{CONTACT M.~W. Hess Email: mhess@sissa.it, A. Quaini Email: quaini@math.uh.edu, G. Rozza Email: gianluigi.rozza@sissa.it} and Annalisa Quaini\textsuperscript{b} and Gianluigi Rozza\textsuperscript{a}}
\affil{\textsuperscript{a}SISSA Mathematics Area, mathLab, International School for Advanced Studies, via Bonomea 265, I-34136 Trieste, Italy; 
\textsuperscript{b}Department of Mathematics, University of Houston, Houston, Texas 77204, USA}
}

\maketitle

\begin{abstract}
We consider the \emph{Navier-Stokes} equations in a channel with a narrowing and walls of varying curvature. By applying the \emph{empirical interpolation method} to generate an affine parameter dependency, the offline-online procedure can be used to compute reduced order solutions for parameter variations. The reduced order space is computed from the steady-state snapshot solutions by a standard POD procedure. The model is discretised with high-order spectral element \emph{ansatz} functions, resulting in 4752 degrees of freedom. The proposed reduced order model produces accurate approximations of steady-state  solutions for a wide range of geometries and kinematic viscosity values. The application that motivated the present study is the onset of asymmetries (i.e., symmetry breaking bifurcation) in blood flow through a regurgitant mitral valve, depending on the Reynolds number and the valve shape. Through our computational study, we found that the critical Reynolds number for the symmetry breaking increases as the wall curvature increases.
\end{abstract}


\section{Introduction and Motivation}

\begin{figure}
\begin{center}
\includegraphics[height = 0.25\textwidth]{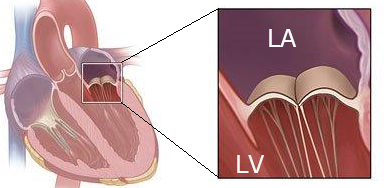}\hskip .7cm
 \includegraphics[height = 0.25\textwidth]{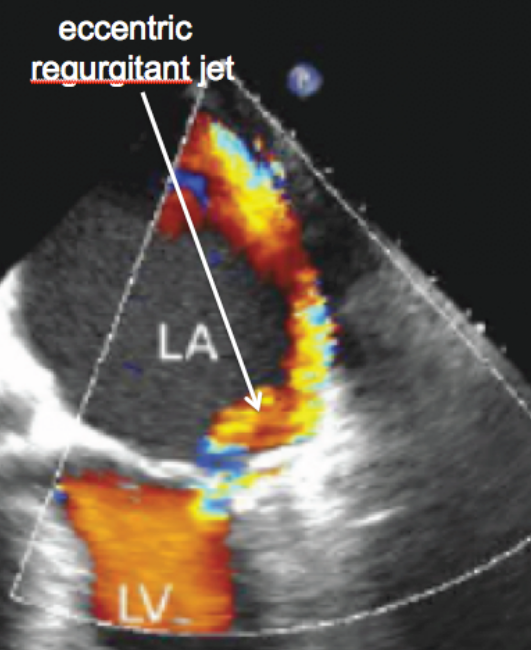}
\end{center}
        \caption{Left: Anatomy of the heart showing the mitral valve.
Right: Echocardiographic image of a jet
flowing from the left ventricle (LV) to the left atrium (LA)
hugging the walls of the LA. Colors denote different fluid velocities.}\label{MitralValve}
\end{figure}

We consider the flow of an incompressible fluid through a planar channel with a narrowing, 
where the walls creating the narrowing have variable curvature.
An application that motivated the present study is the flow of blood through a
regurgitant mitral valve.
Mitral regurgitation (MR) is a valvular disease characterized by abnormal leaking of blood
through the mitral valve from the left ventricle into the left atrium of the heart. See Fig.~\ref{MitralValve}.
In certain cases,  the regurgitant jet ``hugs" the wall of the heart's atrium as shown in  Fig.~\ref{MitralValve} (right).
These wall-hugging, non-symmetric regurgitant jets have been observed at low
Reynolds numbers (\cite{regurgRe1000,regurgRe50}) and are said to undergo the Coanda effect
(\cite{Hess:wille_fernholz_1965}). Such jets represent one of the biggest challenges in echocardiographic 
assessment of MR (\cite{Hess:Ginghina2007}).
In (\cite{Hess:Quaini2016,Hess:Pitton2017534,Hess:Pitton2017,Hess:Wang2017,Hess:max_bif,Hess:ICOSAHOM2018}), 
we made a connection between the cardiovascular and bioengineering 
literature reporting on the Coanda effect in MR
and the fluid dynamics literature with the goal of
identifying and understanding the main features of the corresponding flow conditions.

In (\cite{Hess:Quaini2016,Hess:Pitton2017534,Hess:Pitton2017,Hess:Wang2017,Hess:max_bif,Hess:ICOSAHOM2018}), 
we studied what triggers the Coanda effect in a simplified setting by reformulating
the problem in terms of the hydrodynamic stability of solutions of the incompressible 
Navier-Stokes equations in contraction-expansion channels with straight walls. Such channels 
have the same geometric features of MR
and the wall-hugging effect is nothing but a symmetry breaking bifurcation.
Here, we extend the study to contraction-expansion channels with curved walls
as a first step towards more realistic geometries and eventually fluid-structure interaction. 

For numerical treatment of the incompressible Navier--Stokes problem, we apply 
the spectral element method (SEM), which uses high-order polynomial \emph{ansatz} functions such as Legendre polynomials. See, e.g., (\cite{Hess:PateraSEM, Hess:CHQZ1, Hess:CHQZ2}), 
and (\cite{Hess:Sherwin:2005, Hess:Taddei}) for applications in fluid dynamics.
With a coarse partitioning of the computational domain into spectral elements, the high-order \emph{ansatz} functions are prescribed over each element.
The \emph{ansatz} functions are modified for numerical stability and to enable continuity across element boundaries.
Let  $p$ be the order of the polynomial.
Typically, an exponential error decay under $p$-refinement
can be observed, which provides computational advantages over more standard finite element methods.

Varying wall curvature and kinematic viscosity are considered for the parametric model order reduction.
From a set of sampled high-order solves, a reduced order model is generated, which approximates the high-order solutions and
allows fast parameter sweeps of the two-dimensional parameter domain. 
The offline-online decomposition required for fast reduced order parameter evaluations is established with the \emph{empirical interpolation method}
(\cite{Hess:Barrault}, \cite{Hess:Chaturantabut}, \cite{Hess:Quarteroni2007}, \cite{Hess:Rozza2009}, \cite{Hess:Maday2015}).
The reduced-order modeling (ROM) techniques described in this work are implemented 
in open-source project ITHACA-SEM\footnote{\texttt{https://github.com/mathLab/ITHACA-SEM} and \texttt{https://mathlab.sissa.it/ITHACA-SEM}}.
This extends our previous work (\cite{Hess:ICOSAHOM2018}, \cite{Hess:ENUMATH17_me}) to bifurcations in geometries with non-affine geometry variations.

The outline of the paper is as follows. In sec.~2, the model problem is defined and the parametric variations are explained. 
Sec.~3 provides details on the spectral element discretization and sec.~4 explains the model order reduction with the empirical interpolation.
Numerical results are presented in sec.~5, while in sec.~6 conclusions are drawn and future perspectives and developments are discussed.

\section{Problem Formulation}

Let $\Omega \in \mathbb{R}^2$ be the computational domain.
Incompressible, viscous fluid motion in spatial domain $\Omega$ over a time interval $(0, T)$
is governed by the incompressible
 \emph{Navier-Stokes} equations:
\begin{eqnarray}
\frac{\partial \mathbf{u}}{\partial t} + \mathbf{u} \cdot \nabla \mathbf{u} &=& - \nabla p + \nu \Delta \mathbf{u} + \mathbf{f}, \label{Hess:NSE0} \\
\nabla \cdot \mathbf{u} &=& 0,
\label{Hess:NSE1}
\end{eqnarray}
where $\mathbf{u}$ is the vector-valued velocity, $p$ is the scalar-valued
pressure, $\nu$ is the kinematic viscosity and $\mathbf{f}$ is a body forcing.
\noindent Boundary and initial conditions are prescribed as 
\begin{eqnarray}
\mathbf{u} &=& \mathbf{d} \quad \text{ on } \Gamma_D \times (0, T), \\
\nabla \mathbf{u} \cdot \mathbf{n} &=& \mathbf{g} \quad \text{ on } \Gamma_N \times (0, T), \\
\mathbf{u} &=& \mathbf{u}_0 \quad \text{ in } \Omega \times 0,
\label{Hess:NSE_boundaryCond}
\end{eqnarray}

\noindent with $\mathbf{d}$, $\mathbf{g}$ and $\mathbf{u}_0$ given and $\partial \Omega = \Gamma_D \cup \Gamma_N$, $\Gamma_D \cap \Gamma_N = \emptyset$.
The \emph{Reynolds} number $Re$, which characterizes the flow (\cite{Hess:Holmes}), depends on $\nu$,
a characteristic velocity $U$, and a characteristic length $L$:
\begin{equation}\label{eq:re}
Re = \frac{UL}{\nu}.
\end{equation}

We are interested in the steady states, i.e., solutions where $\frac{\partial \mathbf{u}}{\partial t}$ vanishes.
The high-order simulations are obtained through time-advancement, while the reduced order solutions are computed 
through fixed-point iterations.

\subsection{Non-linear solver}

The \emph{Oseen}-iteration is a secant modulus fixed-point iteration, which in general exhibits a linear rate of convergence (\cite{Hess:Oseen}).
It solves for a steady-state solution, i.e., $\frac{\partial \mathbf{u}}{\partial t} = 0$ is assumed.
Given a current iterate (or initial condition) $\mathbf{u}^k$, the next iterate $\mathbf{u}^{k+1}$ is found
by solving the following linear system:
\begin{eqnarray}
 -\nu \Delta \mathbf{u}^{k+1} + (\mathbf{u}^k \cdot \nabla) \mathbf{u}^{k+1} + \nabla p &=& \mathbf{f}  \text{ in } \Omega, \label{Hess:eq_Oseen_main} \cl
\nabla \cdot \mathbf{u}^{k+1} &=& 0   \text{ in } \Omega, \cl
\mathbf{u}^{k+1} &=& \mathbf{d}  \text{ on } \Gamma_D, \cl
\nabla \mathbf{u}^{k+1} \cdot \mathbf{n} &=& \mathbf{g} \text{ on } \Gamma_N. \el
\end{eqnarray}
Iterations are stopped when the relative difference between iterates falls below a predefined tolerance in a suitable norm, like the $L^2(\Omega)$ or $H^1_0(\Omega)$ norm. 
  
\subsection{Model Description}

\begin{figure}[ht]
\begin{center}
 \includegraphics[scale=.19]{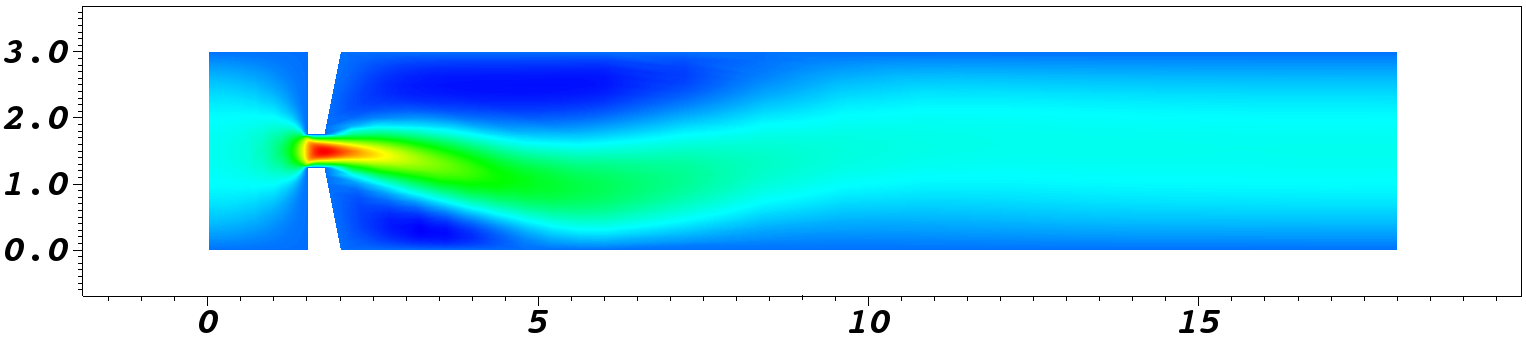} $\quad$
 \includegraphics[scale=.25]{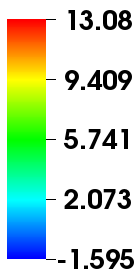} \\
\hspace{.2cm} \includegraphics[scale=.19]{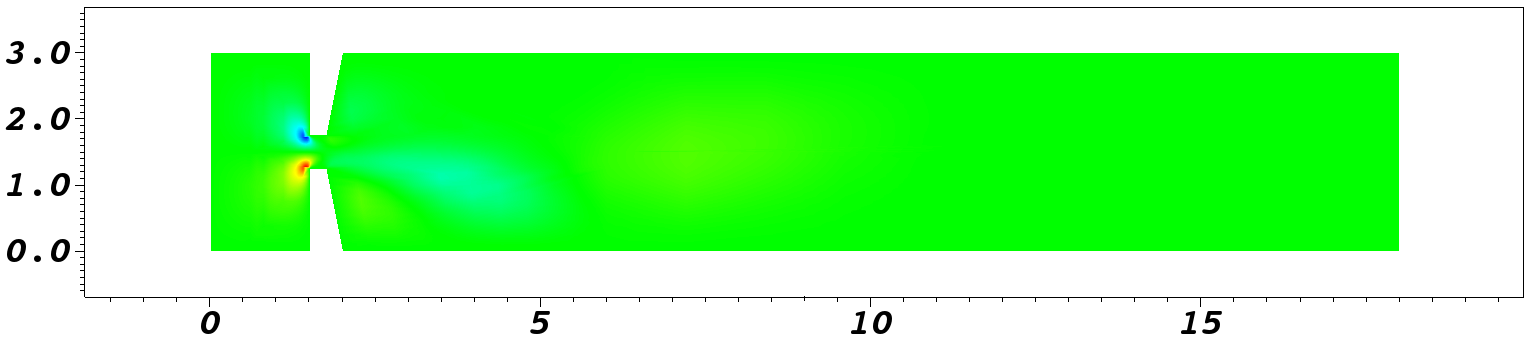} $\quad$
 \includegraphics[scale=.25]{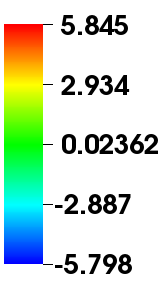}
 \end{center}
 \caption{Full order, steady-state solution in the geometry with straight walls and for $\nu = 0.15$: velocity in x-direction (top) and y-direction (bottom).}
 \label{Hess:FOM_straight}
\end{figure}

\begin{figure}[ht]
\begin{center}
 \includegraphics[scale=.19]{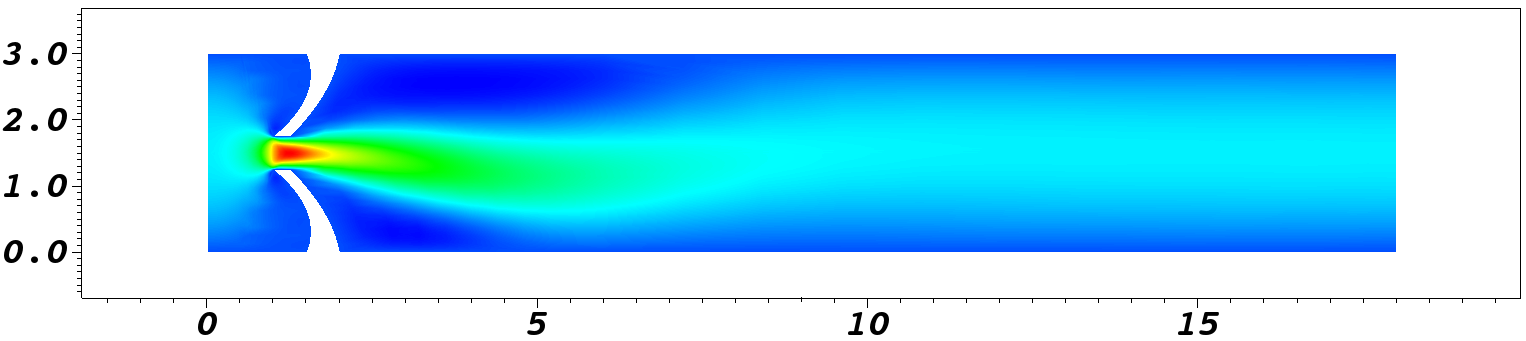} $\quad$
 \includegraphics[scale=.25]{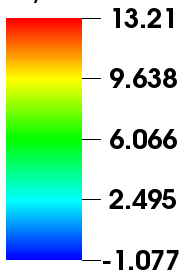} \\
 \hspace{.2cm} \includegraphics[scale=.19]{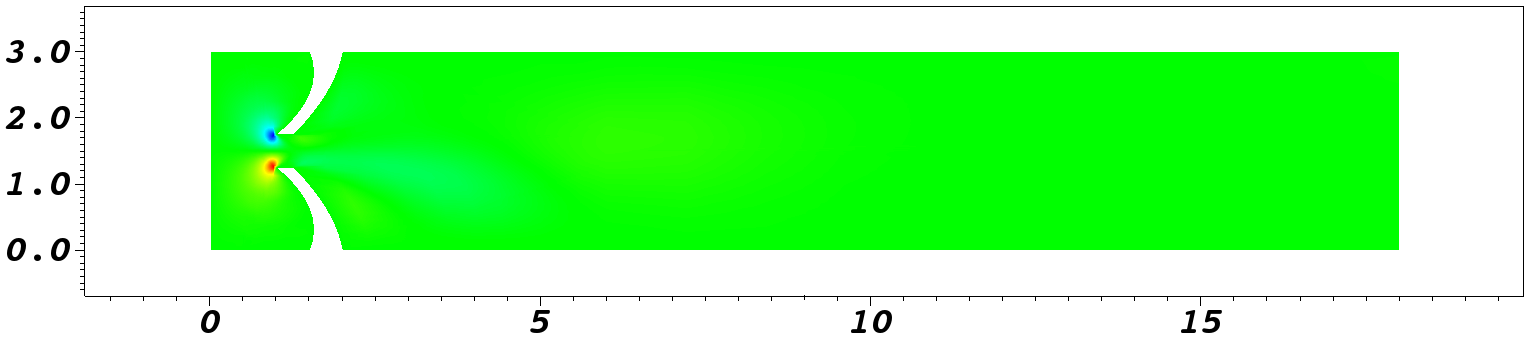} $\quad$
 \includegraphics[scale=.25]{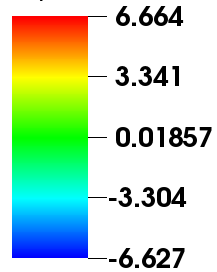}
  \end{center}
 \caption{Full order, steady-state solution in the geometry with curved walls with the largest considered curvature and for $\nu = 0.15$ : velocity in x-direction (top) and y-direction (bottom).}
 \label{Hess:FOM_curved}
\end{figure}

We consider the channel flow through a narrowing created by walls of varying curvature and with variable kinematic viscosity.
See Fig.~\ref{Hess:FOM_straight} and \ref{Hess:FOM_curved} for the steady-state velocity components
for $\nu = 0.15$ in a geometry with straight walls and curved walls, respectively.
In all the cases under consideration, the spectral element expansion uses modal Legendre polynomials of order $p = 10$ for the velocity. 
The pressure \emph{ansatz} space is chosen of order $p-2$ to fulfill the inf-sup stability condition (\cite{Hess:infsup,Hess:BBF}).
A parabolic inflow profile is  prescribed at the inlet (i.e., $x = 0$) with horizontal velocity component 
$u_x(0,y) = y(3-y)$ for $y \in [0, 3]$.
At the outlet (i.e., $x = 18$) we impose a stress-free boundary condition, 
while everywhere else a no-slip condition is prescribed.
We consider symmetric boundary conditions, because we want to study
the symmetry breaking due to the nonlinearity in problem \eqref{Hess:NSE0} - \eqref{Hess:NSE1}.
For a more realistic setting one would have to account for 
different inlet velocity profiles and the pulsatility of the flow 
(i.e., include the Strouhal number among the parameters).


Each curved wall is defined by a second order polynomial, interpolating three prescribed points. 
While the points at the domain boundary $y=0$ and $y=3$ are kept fixed, 
the inner points are moved towards $x=0$ in order to create an increasing curvature.
The viscosity varies in the interval $\nu \in [0.15, 0.2]$.
We recall that the Reynolds number $Re$ \eqref{eq:re} depends on the kinematic viscosity. 
As $Re$ is varied for each fixed geometry, a supercritical pitchfork bifurcation occurs:
for $Re$ higher than the critical bifurcation point, three solutions exit.
Two of these solutions are stable, one with a jet towards the top wall and one with a jet
towards the bottom wall, and one is unstable. The unstable solution 
is symmetric to the horizontal centerline at $y=1.5$, while the jet of the stable solutions
is said to undergo the Coanda effect.

In this investigation, we do not deal with recovering all bifurcation branches, but limit our attention 
to the stable branch of solutions with jets hugging the bottom wall. However, we remark that
recovering all bifurcating solutions with model reduction methods is also possible. See, e.g., (\cite{Hess:Herrero2013132}).

\section{Spectral Element Full Order Discretization}

The spectral/hp element software framework we use for the numerical solution of problem \eqref{Hess:NSE1}
is Nektar++, version 4.4.0\footnote{See \textbf{www.nektar.info}.}.
The large-scale discretized system that has to be solved at each step of the \emph{Oseen}-iteration can be written as
\begin{eqnarray}
\begin{bmatrix}
\begin{array}{ccc}
A & -D^T_{bnd} & B  \\
-D_{bnd} & 0 & -D_{int} \\
\tilde{B}^T & -D^T_{int} & C
\end{array}
\end{bmatrix}
\begin{bmatrix}
\begin{array}{ccc}
\mathbf{v}_{bnd} \\
\mathbf{p} \\
\mathbf{v}_{int} 
\end{array}
\end{bmatrix}
&=
\begin{bmatrix}
\begin{array}{ccc}
\mathbf{f}_{bnd} \\
\mathbf{0} \\
\mathbf{f}_{int}
\end{array}
\end{bmatrix},
\label{Hess:fully_expanded}
\end{eqnarray}
\noindent for fixed parameter vector $\boldsymbol{\mu}$, which denotes the geometrical and physical parameters. 
In \eqref{Hess:fully_expanded}, $\mathbf{v}_{bnd}$ and $\mathbf{v}_{int}$ denote the arrays of the  velocity degrees of freedom on the boundary and in the interior of the domain, respectively.
The array of the pressure degrees of freedom is denoted by $\mathbf{p}$. 
The forcing terms on the boundary and interior are $\mathbf{f}_{bnd}$ and $\mathbf{f}_{int}$, respectively.
Next, we explain the matrix blocks.

Matrix $A$ assembles the boundary-boundary velocity coupling, 
$B$ the boundary-interior velocity coupling, 
$\tilde{B}$ the interior-boundary velocity coupling, and
$C$ assembles the interior-interior velocity degree of freedom coupling.
The matrices $D_{bnd}$ and $D_{int}$ assemble the pressure-velocity boundary and pressure-velocity interior contributions.
Due to the varying geometry, each matrix is dependent on the parameter $\boldsymbol{\mu}$.

The linear system \eqref{Hess:fully_expanded} is assembled in local degrees of freedom, i.e., \emph{ansatz} functions with support extending over spectral element 
boundaries are treated seperately for each spectral element. 
See (\cite{Hess:Sherwin:2005}) for detailed explanations.
As a result, matrices $A, B, \tilde{B}, C, D_{bnd}$ and $D_{int}$
have a block structure, with each block corresponding to a spectral element. 
This allows for an efficient matrix assembly since each spectral element is independent from the others, but 
the local degrees of freedom need to be gathered into the global degrees of freedom 
in order to obtain a non-singular system.

The boundary-boundary global element coupling is achieved with the rectangular assembly matrix $M$, which gathers the local boundary degrees of freedom.
Multiplication of the first row of \eqref{Hess:fully_expanded} by $M^T M$ sets the boundary-boundary coupling in local degrees of freedom:
\begin{eqnarray}
\begin{bmatrix}
\begin{array}{ccc}
M^T M A & -M^T M D^T_{bnd} & M^T M B  \\
-D_{bnd} & 0 & -D_{int} \\
\tilde{B}^T & -D^T_{int} & C
\end{array}
\end{bmatrix}
\begin{bmatrix}
\begin{array}{ccc}
\mathbf{v}_{bnd} \\
\mathbf{p} \\
\mathbf{v}_{int} 
\end{array}
\end{bmatrix}
&=
\begin{bmatrix}
\begin{array}{ccc}
M^T M \mathbf{f}_{bnd} \\
\mathbf{0} \\
\mathbf{f}_{int}
\end{array}
\end{bmatrix} .
\label{Hess:fully_expanded_MtM}
\end{eqnarray}

The action of the matrix in \eqref{Hess:fully_expanded_MtM} on the prescribed Dirichlet boundary conditions is computed and added to the source term.
Since the Dirichlet boundary conditions are known, the corresponding equations are removed from the system.
Let $N_\delta$ denote the system size after removal of the known boundary conditions.
The resulting system of high-order dimension $N_\delta \times N_\delta$ is composed of the block matrices and depends on the parameter $\boldsymbol{\mu}$. 
For simplicity of notation, we will write such system in compact form as:
\begin{eqnarray}
\mathcal{A}(\boldsymbol{\mu}) \mathbf{x}(\boldsymbol{\mu}) = \mathbf{f}(\boldsymbol{\mu}).
\label{Hess:final_to_be_proj}
\end{eqnarray}

\section{Reduced Order Space Generation}

The ROM computes an approximation to the full order model using a few modes of the 
POD as \emph{ansatz} functions (\cite{Hess:Lassila2014}).
To achieve a computational speed-up, the matrix assembly for a new parameter of interest is independent of the large-scale discretization size $N_\delta$.
The \emph{empirical interpolation method} (see sec.~\ref{sec:EIM})
computes an affine parameter dependency, which enables an offline-online decomposition (see sec.~\ref{sec:off-on}).
After a time-intensive offline phase, reduced order solves can be evaluated quickly over the parameter range of interest.

The POD of $N$ (typically small) uniformly sampled full-order solves, called \emph{snapshots}, is performed. 
The most dominant modes corresponding to $99.99\%$ of the POD energy (as suggested in \cite{Hess:Lassila2014})
form the projection matrix $U \in \mathbb{R}^{N_\delta \times N}$ and implicitly define the low-order space $V_N = \text{span} (U)$.
The large-scale system \eqref{Hess:final_to_be_proj} is then projected onto the reduced order space:
\begin{eqnarray}
U^T \mathcal{A}(\boldsymbol{\mu}) U \mathbf{x}_N(\boldsymbol{\mu}) = U^T \mathbf{f}(\boldsymbol{\mu}) .
\label{Hess:final_proj}
\end{eqnarray}
The low order solution $\mathbf{x}_N(\boldsymbol{\mu})$ approximates the large-scale solution as $\mathbf{x}(\boldsymbol{\mu}) \approx U \mathbf{x}_N(\boldsymbol{\mu})$.
The stability properties of the full-order model do not necessarily carry over to the reduced-order model, which can introduce instabilities.
In particular, the reduced order inf-sup stability constant might approach zero for some parameter value, while the full-order inf-sup stability constant is bounded 
away from zero. One way to alleviate this problem is by using inf-sup supremizers (\cite{Hess:Lassila2014}) 
or considering space-time variational approaches (\cite{Hess:Yano2013}).

\subsection{Empirical Interpolation Method}\label{sec:EIM}

The discrete \emph{empirical interpolation method} (EIM) (\cite{Hess:Barrault}, \cite{Hess:Chaturantabut}, \cite{Hess:Quarteroni2007}, \cite{Hess:Rozza2009}) 
computes an approximate affine parameter dependency.
During the snapshot computation, the parameter-dependent matrices
are collected. The matrix discrete empirical interpolation (\cite{Hess:Negri}) allows to decompose $\mathcal{A}(\boldsymbol{\mu})$
as follows:
\begin{eqnarray}
 \mathcal{A}(\boldsymbol{\mu}) \approx \sum_{i=1}^{Q_a} \tau_i(\boldsymbol{\mu}) A_i,
\label{Hess:MDEIM}
\end{eqnarray}
where $\tau_i(\boldsymbol{\mu})$ are scalar parameter-dependent coefficient functions and 
$A_i$ are parameter-independent matrices.
Each coefficient function $\tau_i(\boldsymbol{\mu})$ corresponds to a single matrix entry of $\mathcal{A}(\boldsymbol{\mu})$.
Since the assembly of only a few (here $Q_a < 30$) matrix entries can be implemented efficiently, an approximation of $\mathcal{A}(\boldsymbol{\mu})$ is readily available 
for each new $\boldsymbol{\mu}$.
As for the projection space $U$, $99.99\%$ of the POD energy is used to approximate the system matrices from the collected matrices during snapshot computation.
The EIM is preformed for each submatrix identified in \eqref{Hess:fully_expanded_MtM} and the actual system matrix $\mathcal{A}(\boldsymbol{\mu})$ is then composed of the 
separately approximated block matrices.

\subsection{Offline-Online Decomposition}
\label{sec:off-on}

The offline-online decomposition (\cite{Hess:RBref}) enables the computational speed-up of the ROM approach in many-query scenarios.
It relies on an affine parameter dependency, such that all computations depending on the high-order model size 
$N_\delta$ can be performed in a parameter-independent offline phase. 
Then, the input-output evaluation performed online is independent of $N_\delta$ and thus fast.

After applying the \emph{empirical interpolation method} in the geometry parameter, the parameter dependency is cast in an affine form.
Therefore, there exists an affine expansion of the system matrix $\mathcal{A}(\boldsymbol{\mu})$ in the parameter $\boldsymbol{\mu}$ given by \eqref{Hess:MDEIM}.
To achieve fast reduced order solves, the offline-online decomposition
computes the parameter independent projections offline, which are
stored as small-sized matrices of the order $N \times N$.
Since in an \emph{Oseen}-iteration each matrix is dependent on the previous iterate, the submatrices corresponding to each basis function are assembled and then 
formed online using the affine expansion computed from the EIM and a fast evaluation of a single matrix entry as required by the EIM coefficient functions $\tau_i$.

\section{Numerical Results}

Snapshot solutions are sampled over a uniform $8 \times 9$ grid from the full-order model, with $8$ samples along the $\nu$ parameter direction and $9$ samples along the 
geometry parameter direction. 
The number of required snapshot computations might potentially be reduced
when using a greedy sampling, which requires
error indicators or error estimators. See \cite{Hess:RBref}. 
Error estimation does even allow a certification of the ROM accuracy, 
but it requires an estimation of the inf-sup constant as well as a bound on the \emph{empirical interpolation} error.

\begin{figure}[ht]
\centering
 \includegraphics[scale=1]{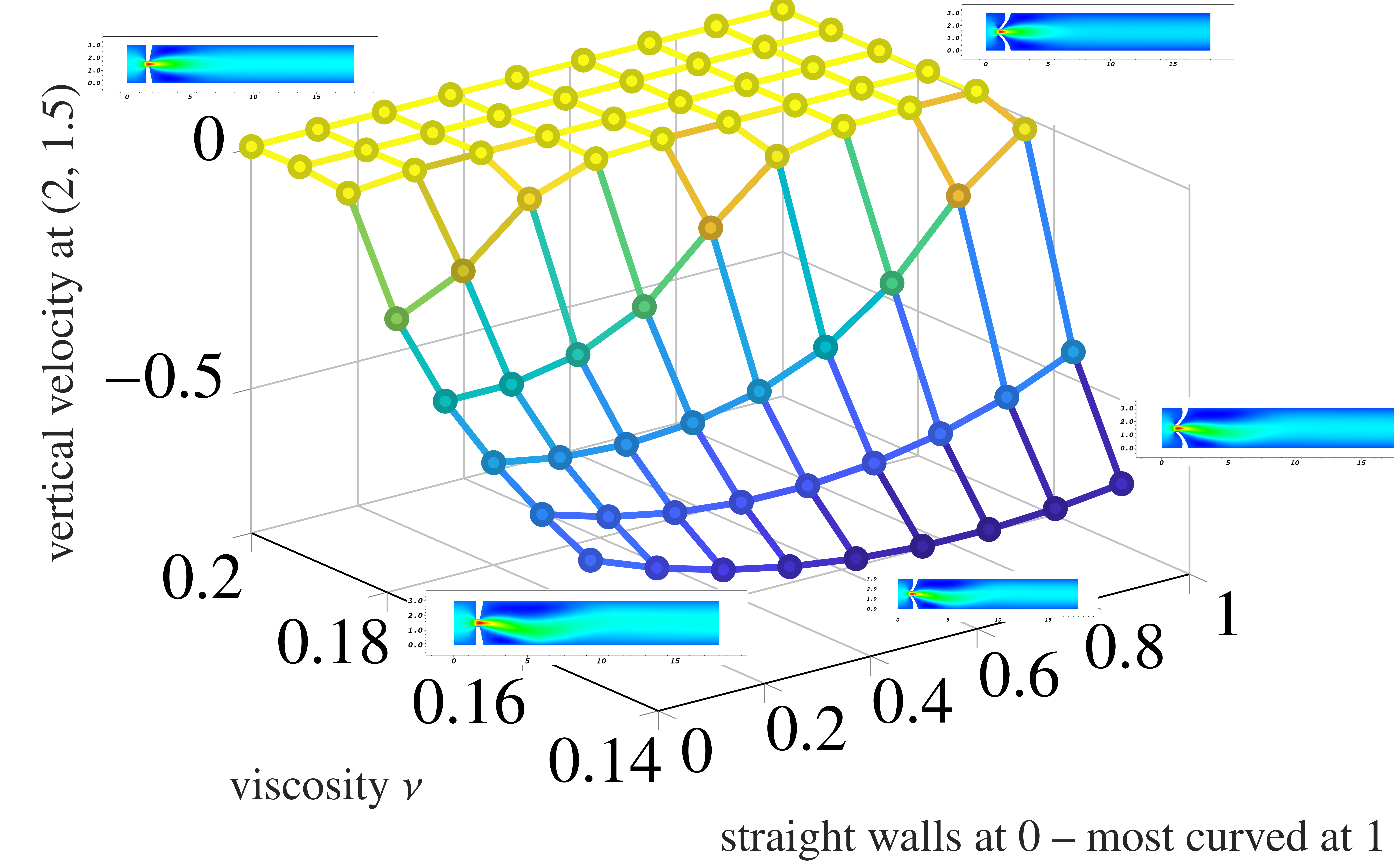}
 \caption{Bifurcation diagram of the full order model over a uniform $8 \times 9$ grid.}
 \label{Hess:FOM_sampled}
\end{figure}

The vertical velocity at the point $(2, 1.5)$ is used to generate the bifurcation diagram
reported in Fig.~\ref{Hess:FOM_sampled}.
As expected, in a fixed geometry the symmetry breaking bifurcation
occurs when the $Re$ exceeds a critical value. It is very interesting to observe that 
such critical value increases 
as the wall curvature increases, i.e.~as the walls get more curved, the stronger the inertial forces need to
be to break the symmetry of the solution.
This means that the estimates for the critical $Re$ in 3D geometries with straight walls
found in \cite{Hess:Pitton2017534,Hess:Wang2017} provide a lower bound for the critical $Re$
at which the Coanda effect is observed in vivo (see Fig.~\ref{MitralValve} (right)).

We remark that the 2D case can be seen as a limit of the 3D case for channel depth tending to 
infinity. In \cite{Hess:Pitton2017534}, the influence of the channel depth on the flow pattern
is investigated. It is shown that non-symmetry jets appear at higher $Re$
as the channel depth gets smaller. Thus, we expect that the critical $Re$ 
for the symmetry breaking in a 3D channel with curved walls to be higher
than the values reported here for a 2D channel.
This indicates the Coanda effect occurs in mitral valves with elongated
orifices (corresponding to deeper channels).


\begin{figure}[ht]
\centering
 \includegraphics[scale=.8]{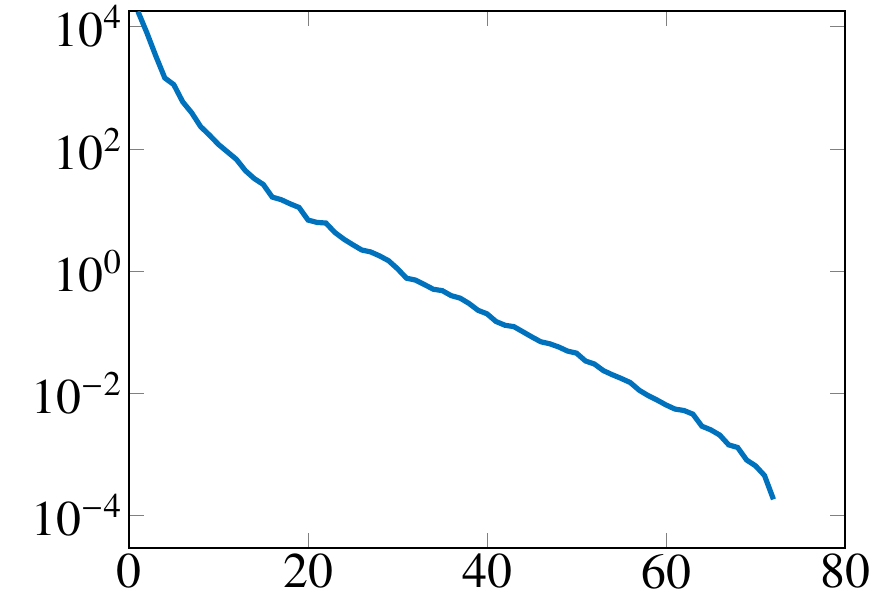}
 \caption{Decay of POD energy of the sampled snapshot solutions.}
 \label{Hess:PODdecay}
\end{figure}

The accuracy of the ROM is assessed using $N = 72$ snapshots for the POD to recover the original snapshot data.
Fig.~\ref{Hess:PODdecay} shows the decay of the energy of the POD modes. 
To reach the typical threshold of $99.99\%$ on the POD energy, $N = 33$ POD modes are required as RB ansatz functions.

\begin{figure}[ht]
\centering
 \includegraphics[scale=.8]{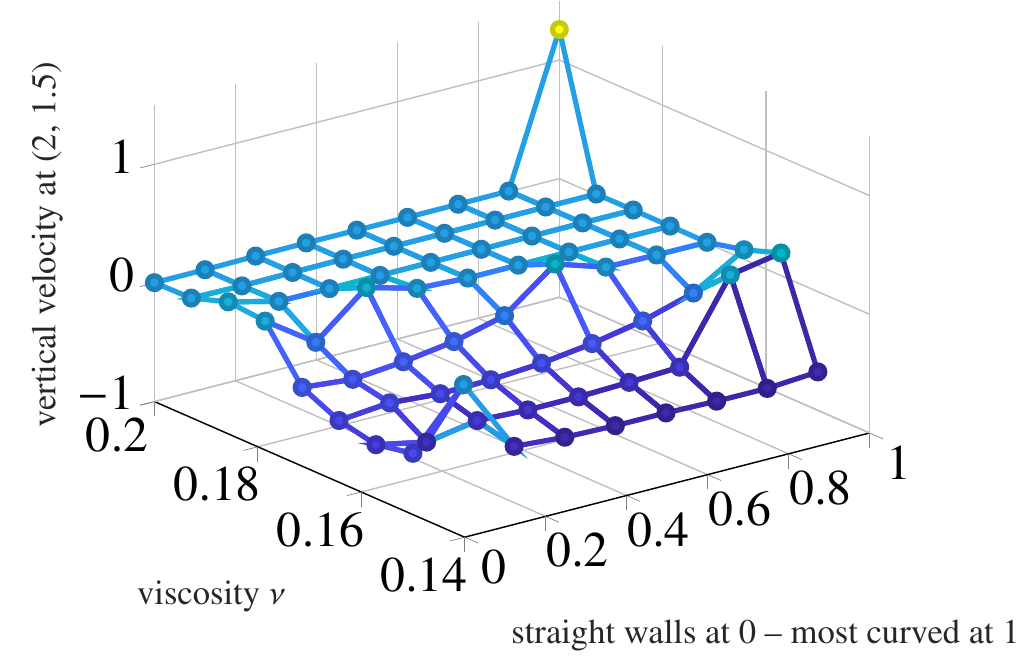}
 \caption{Bifurcation diagram of the reduced order model with reduced model size $20$.}
 \label{Hess:ROM20}
\end{figure}

Fig.~\ref{Hess:ROM20} shows the bifurcation diagram of the reduced order model with $N = 20$ basis functions.
The absolute error at the point value is less than $0.01$ at $46$ parameter locations and less than $0.1$ at $63$ parameter locations.
This indicates, that the high-order solutions have been well-resolved at these configurations.
There are a few outliers, where the iteration scheme did not converge and the value for the last iterate is shown. 
Most likely this can be resolved by taking a finer snapshot sampling into account
or by using a localized reduced-order modeling approach (\cite{Hess:max_bif}).

The \emph{empirical interpolation method} relies on the fast computation of a few matrix entries during the online phase. 
Since the spectral element \emph{ansatz} functions have support over a whole spectral element, this operation cannot be performed as fast as 
with a finite element or finite volume method for instance, where \emph{ansatz} functions have a local support. 
Nevertheless, the computational gain is significant after the affine form has been established.
The time requirement for a single fixed point iteration step reduces from about $10$ s to $0.1$ s.

\section{Conclusion and Outlook}

We proposed a reduced order model that combines 
\emph{empirical interpolation method} and a POD reduced basis technique 
to recover full-order solutions of the Navier--Stokes equations in domains with walls of varying curvature (non-affine variation). 
The non-linear geometry changes allow to simulate more realistic scenarios
in the context of the Coanda effect in cardiology, but also require a fine sampling at the snapshot and \emph{empirical interpolation} level. 
Since the model problem studied here undergoes a supercritical pitchfork bifurcation, 
introducing non-unique solutions, further numerical techniques are required to recover all bifurcation branches.
The spectral element method is a suitable method for these tasks.  
However, the computational gain that one can expect is not as significant
as in the case of methods using \emph{ansatz} functions with a local support, such as 
the finite element method.

As a next step, we will enhance the reduced order model proposed here
by using localizes bases in order to recover every solutions in the considered parameter domain
with high accuracy.

\section*{Acknowledgments}  

This work was supported by European Union Funding for Research and Innovation through the European Research Council
(project H2020 ERC CoG 2015 AROMA-CFD project 681447, P.I. Prof. G. Rozza).
This work was also partially supported by NSF through grant DMS-1620384 (A. Quaini).

\end{document}